\newcommand{\be}{\begin{eqnarray}}
\newcommand{\ee}{\end{eqnarray}}

\def\fr{\frac{1}{2}}
\def\mref#1{(\ref{#1})}
\def\zero{\setcounter{equation}{0}}

\def\p{\partial}

\def\bd{\begin{displaymath}}
\def\ed{\end{displaymath}}

\def\ba#1{\begin{array}{#1}}
\def\ea{\end{array}}
\def\nn{\nonumber}
\newfont{\Bbb}{msbm10 scaled 1200}
\documentclass[12pt]{article}

\title{
Poisson--Lie structures on Galilei group\footnote{supported by \L\'od\'z University grant no 698}
}
\author{
{Y.Brihaye}\\
{Department of Mathematical Physics, University of Mons,} \\
{Av.Maistriau, B-7000 Mons, Belgium.}\\
{}\\
{}\\
{E.Kowalczyk, P.Ma\'slanka}\\
{Department of Theoretical Physics II, University of \L\'od\'z,}\\
{ul. Pomorska 149/153, 90236 \L\'od\'z, Poland.}
}
\date{}
\begin{document}
\maketitle{}
\thispagestyle{empty}
\begin{abstract}
The complete list of Poisson--Lie structures on 4-d Galilei group is presented.
\end{abstract}
\newpage
\section{Introduction \label{s1}}

\indent Quantum groups have emerged in physics in connection with an attempt to
understand the symmetries underlying exact solvability of certain 
quantum--mechanical and statistical models; they appeared to be quite powerful
in this respect.

Therefore, it is natural to ask whether their range of applicability as 
a~mathematical tool for describing physical symmetries is wider and covers,
in particular, the most important case of space--time symmetries. The 
theory of Hopf algebras offers a variety of structures which can be viewed
as deformations of classical space--time symmetry groups. For example,
a number of deformed Poincare groups were considered~\cite{b1}.
They posses many attractive features. However, if one is going to take 
seriously the very idea of quantum space time symmetries, many conceptual
problems arise which solution seems to be quite difficult. It is rather obvious 
that one should concentrate on deformations of relativistic symmetries because
it is  a high energy/ small distance region where the deviations from
the predictions of ``classical'' theory should occur. However, as far as the 
conceptual problems are concerned, the nonrelativistic region provides similar 
challenge. On the other hand it seems that the study of nonrelativistic deformed
symmetries is slightly simpler. One is not here faced with some complications
typical for relativistic quantum theory as, for example, the nonexistence
of particle number conserving nontrivial interactions.

In the present paper as a first step towards the understanding of 
nonrelativistic quantum space--time symmetries we classify all 
nonequivalent Poisson--Lie structures on 
the four-dimensional Galilei group.
The method we use (c.f.\cite{b2}) is based on solving directly the cocycle
condition; it has been already used for finding Poisson--Lie structures
on the two-dimensional Galilei group~\cite{b3},\cite{b4}.

We find families of structures which cannot be related to each other by the
automorphisms of the Galilei group; contrary to the case of the
Poincare group \cite{b5}
many of them are not of the coboundary type. 

The paper is organized as follows. 
In sec.\ref{s2}\ we sketch a~general strategy while the results are presented
in sec.\ref{s3}. All technical details are relegated to the number of appendices.

Let us conclude the introduction which some details concerning 
the Galilei\\ group\cite {b6}.

The generic element $g$\ of the 
ten--parameter Galilei group $G$\ is denoted by
\be
g&=&(t,\vec{a},\vec{v},R)\label{w1}
\ee
We shall denote by\\
$T=\{(t,\vec{0},\vec{0},I)\}$\\
$S=\{(0,\vec{a},\vec{0},I)\}$\\
$V=\{(0,\vec{0},\vec{v},I)\}$\\
${\cal R}=\{(0,\vec{0},\vec{0},R)\}$\\
the subgroups of time translations, space--translations, pure Galilei transformations (boosts),
rotations, respectivly.

The group law is expressed by
\be
g''=g'\cdot g=(t',\vec{a'},\vec{v'},R')(t,\vec{a},\vec{v},R)=\nn\\
=(t+t',\vec{a'}+R'\vec{a}+\vec{v'}t,\vec{v'}+R'\vec{v},R'R)\label{w2}
\ee
The identity for the group is
\be
e=(0,\vec{0},\vec{0},I)\label{w3}
\ee
and  the inverse of the generic element is given by
\be
g^{-1}=&(t,\vec{a},\vec{v},R)^{-1}&=(-t,-R^{-1}(\vec{a}-t\vec{v}),-R^{-1}\vec{v},R^{-1})\label{w4}
\ee
The generators $H,\vec{P},\vec{K}$\ and $\vec{J}$\ of 
the Galilei Lie algebra
${\cal G}$\ are defined with the help of the exponential parametrization
\be
g&=&e^{-itH}e^{i\vec{a}\vec{P}}e^{i\vec{v}\vec{K}}e^{i\vec{\theta}\vec{J}}\label{w5}
\ee
and they obey the
following commutation rules (only the nonvanishing ones are written up)
\be
{}[J_i,J_j]&=&i\varepsilon_{ijk}J_k\nn\\
{}[J_i,K_j]&=&i\varepsilon_{ijk}K_k\label{w6}\\
{}[J_i,P_j]&=&i\varepsilon_{ijk}P_k\nn\\
{}[K_i,H]&=&iP_i\nn\\
\ee
here and in the sequel the summation over repeated indices is understood
$(i,j,k=1,2,3)$.

The automorphism group of $G$\ consists of the inner automorphisms together
with two outer ones generated by space dilations and time dilations~:
\be
(t,\vec{a},\vec{v},R)&\to&(t,a\vec{a},a\vec{v},R)\nn\\
(t,\vec{a},\vec{v},R)&\to&(bt,\vec{a},b^{-1}\vec{v},R)\label{w7}
\ee
In what follows we shall need the right-invariant vector fields on $G$. 
Denoting by $\hat{X}$\ the right-invariant vector field corresponding to the
element $X$\ of the Lie algebra ${\cal G}$\ we have
\be
\hat{H}&=&-i{\p\over\p t}\nn\\
\hat{P_i}&=&i{\p\over\p a_i}\label{w8}\\
K_i&=&i(t{\p\over\p a_i}+{\p\over\p v_i})\nn\\
J_i&=&-i\varepsilon_{ijk}a_j{\p\over\p a_k}-i\varepsilon_{ijk}v_j{\p\over\p v_k}
-i\varepsilon_{ijk}R_{jl}{\p\over\p R_{kl}}\nn
\ee
\section{Poisson--Lie structures on Galilei group --- the general strategy\label{s2}}

\indent Let us remind the notion of Poisson--Lie group~\cite{b7},\cite{b8}.
It is a Lie group $\tilde{G}$\ which has a~Poisson structure \hbox{\{ , \}}\ 
such that the multiplication map \hbox{$m: \tilde{G}\times\tilde{G}\to\tilde{G}$}\ 
is a~Poisson map (here \hbox{$\tilde{G}\times\tilde{G}$}\ is given 
by the product
 Poisson structure).

Poisson--Lie structures can be described explicitly as follows. Let ${\tilde{\cal G}}$\ 
be the Lie algebra of $\tilde{G}$; denote by $\{X_i\}$\ an arbitrary basis in 
${\tilde{\cal G}}$\ and  let $c_{ij}^k$\ be the corresponding structure constants.
One defines a~mapping \hbox{$\eta:\tilde{G}\to\wedge^2\tilde{\cal G}$}, 
\be
\eta(g)\equiv\eta^{ij}(g)X_i\otimes X_j &,& \eta^{ij}(g)=-\eta^{ji}(g)\label{w9}
\ee
Let $\{X^R_i\}$\ be the realization of $\tilde{\cal G}$\ in terms of right
invariant vector fields on~$\tilde{G}$. The Poisson bracket on $\tilde{G}$\ 
given by
\be
\{\Phi,\Psi\}&=&\eta^{ij}(X^R_i\Phi)(X_j^R\Psi)\label{w10}
\ee
defines the Poisson--Lie structure on $\tilde{G}$\ provided the following conditions are
obeyed\\
(i) Poisson--Lie property (co-cycle condition)
\be
\eta(g'g)&=&\eta(g')+Ad(g')\eta(g)\label{w11}
\ee
(ii) Jacobi identity
\be
\eta^{il}X^R_l\eta^{jk}+\eta^{kl}X^R_l\eta^{ij}+\eta^{jl}X^R_l\eta^{ki}-\nn\\
{}-c_{lp}^j\eta^{il}\eta^{pk}-c^i_{lp}\eta^{kl}\eta^{pj}-
c^k_{lp}\eta^{jl}\eta^{pi}=0\label{w12}
\ee

The inverse is also true: any Poisson--Lie structure on $\tilde{G}$\ can be 
written in the above form. The infinitesimal analogues of Poisson--Lie
groups are Lie bialgebras. For any \hbox{$X\in\tilde{\cal G}$} define
\be
\delta(X)&\equiv&{d\over dt}\eta(e^{itX})|_{t=0}\label{w13}
\ee
Then it can be easily shown that $\delta:\tilde{\cal G}\to \wedge^2\tilde{\cal G}$\ 
has the following properties which are the infinitesimal counterparts of~(i) and~(ii):\\
(i') co-cycle conditions
\be
\delta([X,Y])&=&[X\otimes I+I\otimes X,\delta(Y)]+[\delta(X),I\otimes Y+Y\otimes I]\label{w14}
\ee
(ii') co-Jacobi identity
\be
\sum\limits_{c.p.}(\delta\otimes id)\circ\delta(X)&=&0\label{w15}
\ee
where $c.p.$\ means the summation over cyclic permutation of the factors in
\hbox{$\tilde{\cal G}\otimes\tilde{\cal G}\otimes\tilde{\cal G}$}.

Every Poisson--Lie structure on $\tilde{G}$\ defines the bialgebra structure
 on~$\tilde{\cal G}$. The inverse is true provided $\tilde{G}$\ is connected and
simply connected\cite{b9}. Two Poisson--Lie structures on~$\tilde{G}$\ will
be called equivalent if there exists an automorphism of~$\tilde{G}$\ which is a 
Poisson map.

The main aim of the present paper is to classify, up to equivalence, all
Poisson--Lie structures on the four-dimensional 
Galilei group~$G$. We adopt the
following, rather straightforward, strategy.

First, we write out $\eta$\ in the form
\be
\eta(g)&=&\Psi_i(g)H\wedge J_i+\Phi_i(g)H\wedge P_i+\Gamma_i(g) H\wedge K_i\nn\\
&&+\Lambda_i(g)\varepsilon_{ijk}P_j\wedge P_k+\Upsilon_{ij}(g)P_i\wedge K_j\label{w16}\\
&&+\Sigma_{ij}(g)P_i\wedge J_j+\Xi_i(g)\varepsilon_{ijk}K_j\wedge K_k+\nn\\
&&+\Omega_{ij}(g)K_i\wedge J_j+\Pi_i(g)\varepsilon_{ijk}J_j\wedge J_k\nn
\ee
where 
\hbox{$g\equiv (t,\vec{a},\vec{v},R)$} is an arbitrary element of 
the Galilei group~$G$.

Inserting the expansion \mref{w16} into the one-cocycle condition \mref{w11} one
obtains the following set of functional equations for the coefficients $\Psi_i,\Phi_i$\ , etc.
\be
\Psi_i(gg')&=&\Psi_i(g)+R_{il}\Psi_l(g')\nn\\
\Phi_i(gg')&=&\Phi_i(g)+R_{il}\Phi_l(g')+\varepsilon_{ink}(v_nt-a_n)R_{kl}\Psi_l(g')-\nn\\
&&-tR_{il}\Gamma_l(g')\nn\\
\Gamma_i(gg')&=&\Gamma_i(g)+R_{il}\Gamma_l(g')-\varepsilon_{ink}v_nR_{kl}\Psi_l(g')\nn\\
\Lambda_i(gg')&=&\Lambda_i(g)+R_{il}\Lambda_l(g')-{1\over 2}\varepsilon_{imn}v_mR_{nl}\Phi_l(g')+\nn\\
&&+\fr [(\vec{v}^2t-a_mv_m)\delta_{in}-v_n(v_it-a_i)]R_{nl}\Psi_l(g')+\nn\\
&&+\fr t\varepsilon_{imn}v_mR_{nl}\Gamma_l(g')-\fr t R_{il}\varepsilon_{lmn}\Upsilon_{mn}(g')+\nn\\
&&+\fr [(tv_i-a_i)\delta_{mp}-R_{im}R_{np}(tv_n-a_n)](\Sigma_{pm}(g')\nn\\
&&-\fr t\Omega_{pm}(g'))+t^2R_{il}\Xi_l(g')+\nn\\
&&+(t^2v_pv_i-v_pa_it-v_ia_pt+a_pa_i)R_{pm}\Pi_m(g')\nn\\
\Upsilon_{ij}(gg')&=&\Upsilon_{ij}(g)+R_{im}R_{jn}\Upsilon_{mn}(g')+v_i\varepsilon_{jnk}v_nR_{kl}\Psi_l(g')\nn\\
&&-v_iR_{jl}\Gamma_l(g')-\varepsilon_{jnl}v_nR_{ip}R_{lk}\Sigma_{pk}\nn\\
&&-2t\varepsilon_{ijs}R_{sl}\Xi_l(g')\nn\\
&&+[\varepsilon_{inl}(a_n-v_nt)R_{jp}+\varepsilon_{jnl}v_ntR_{ip}]R_{lk}\Omega_{pk}(g')\nn\\
&&+2[\varepsilon_{ijn}v_n(a_p-v_pt)R_{pm}-\varepsilon_{njr}a_nv_rR_{im}]\Pi_m(g')\label{w17}\\
\Sigma_{ij}(gg')&=&\Sigma_{ij}(g)+R_{ip}R_{jk}\Sigma_{pk}(g')-v_iR_{jl}\Psi_l(g')\nn\\
&&-tR_{ip}R_{jk}\Omega_{pk}(g')+\nn\\
&&+2[(a_j-tv_j)R_{im}-(a_l-tv_l)R_{lm}\delta_{ij}]\Pi_m(g')\nn\\
\Xi_i(gg')&=&\Xi_i(g)+R_{in}\Xi_n(g')-\fr (v_i\delta_{pm}-R_{im}R_{np}v_n)\Omega_{pm}(g')+\nn\\
&&+v_iv_pR_{pm}\Pi_m(g')\nn\\
\Omega_{ij}(gg')&=&\Omega_{ij}(g)+R_{ip}R_{jk}\Omega_{pk}(g')+2(R_{im}v_j-R_{lm}v_l\delta_{ij})\Pi_m(g')\nn\\
\Pi_i(gg')&=&\Pi_i(g)+R_{im}\Pi_m(g')\nn
\ee

In spite of their complicated structure they can be solved in the following way
(cf.\cite{b2},\cite{b3},\cite{b4}). One decomposes the general element 
\hbox{$g=(t,\vec{a},\vec{v},R)$} into the product of four elements belonging
to the subgroups of time--and--space--translations, boosts and rotations, 
(see sec.\ref{s1}).
\be
(t,\vec{a},\vec{v},R)&=&(t,\vec{0},\vec{0},I)\cdot(0,\vec{a},\vec{0},I)\cdot(0,\vec{0},\vec{v},I)\cdot
(0,\vec{0},\vec{0},R)\label{w18}
\ee
According to the condition (i) one can successively calculate $\eta(g)$\ using the
above decomposition provided the form of $\eta$\ for all four subgroups is 
known. In order to find the latter we specify eqs.\mref{w17} to those 
subgroups. The resulting equations can be easily solved; this is done in Appendix~A.
However, in obtaining the final form of $\eta$\ we apply eq.(19) with some definite 
order of multiplication
(for example using \hbox{$(t,\vec{a},\vec{v},R)=(t,\vec{0},\vec{0},I)\cdot$}\\
\hbox{$\cdot((0,\vec{a},0,I)((0,0,\vec{v},I)(0,\vec{0},\vec{0},R)))$})
so there could be further constraints on the
parameters entering $\eta$\ following 
from associativity. Therefore, we reinsert $\eta$\ into eq.\mref{w11} with arbitrary
$g$\ and $g'$\ to find all missed relations between parameters. In this way we 
produce the general solution to eq.\mref{w11} described in Appendix~B (see eq.\mref{wb1}).

There remains to solve (ii) which imposes additional relations between 
coefficients of $\eta$. It is very tedious to try to solve eq.\mref{w12} directly
so we adopt a~different method. From our general form of $\eta$\ we calculate
$\delta$\ and impose (ii') which, in this context, is equivalent to (ii). On 
the other hand it is well known that (ii') can be restated as the condition 
that the dual map $\delta^*$\ defines a Lie algebra structure on ${\cal G}^*$.
Therefore, we first calculate $\delta$\ and the commutators on ${\cal G}^*$\ 
resulting from it and then we solve the Jacobi identities. This is still 
a complicated problem but it can be simplified by using the boost and translation
automorphisms of the Galilei group/algebra. Once this is done there remains only
to use the residual automorphisms to put our solutions in canonical position.
The more detailed discussion is given in Appendix~C.

At the end, having the form of $\eta$\ (eq.\mref{w16}) and the form of the right-invariant 
vector fields (eq.\mref{w8}), using eq.\mref{w10} and taking into
account that for all families of solutions $\Pi_i(g)\equiv 0$\ (see eq.\mref{wb1}
and sec.\ref{s3}) one can easily calculate the following fundamental 
Poisson--Lie brackets.
\be
\{R_{ab},R_{cd}\}&=&0\nn\\
\{v_a,R_{bc}\}&=&\varepsilon_{bjl}R_{lc}\Omega_{aj}\nn\\
\{a_a,R_{bc}\}&=&\varepsilon_{bjl}R_{lc}(\Sigma_{aj}+t\Omega_{aj})\nn\\
\{t,R_{bc}\}&=&-\varepsilon_{bjl}R_{lc}\Psi_j\nn\\
\{t,v_a\}&=&\Gamma_a-\varepsilon_{ajl}v_l\Psi_j\label{w19}\\
\{v_a,v_b\}&=&-2\varepsilon_{abj}\Xi_j+\varepsilon_{bjl}v_l\Omega_{aj}-\varepsilon_{ajl}v_l\Omega_{bj}\nn\\
\{a_a,v_b\}&=&-\Upsilon_{ab}+\varepsilon_{bjl}v_l\Sigma_{aj}-\nn\\
&&-2t\varepsilon_{abj}\Xi_j+t\varepsilon_{bjl}v_l\Omega_{aj}-\varepsilon_{ajl}a_l\Omega_{bj}\nn\\
\{a_a,a_b\}&=&-2\varepsilon_{abj}\Lambda_j-t\Upsilon_{ab}+t\Upsilon_{ba}-2t^2\varepsilon_{abj}\Xi_j\nn\\
&&+\varepsilon_{bjl}a_l\Sigma_{aj}-\varepsilon_{ajl}a_l\Sigma_{bj}+t\varepsilon_{bjl}a_l\Omega_{aj}\nn\\
&&-t\varepsilon_{ajl}a_l\Omega_{bj}\nn\\
\{t,a_a\}&=& -\varepsilon_{ina}\alpha_n\Psi_i+\Phi_a+\Gamma_a\nn
\ee

\section{Poisson--Lie structures on Galilei group --- the results\label{s3}}

\indent By applying the procedure outlined above we  solve the relevant
Jacobi identities for dual algebra (making use of boost and translation automorphisms) and we
arrive at the following families of Poisson--Lie structures (for all cases
\hbox{$\vec{n}=0$}).
\begin{description}
\item[I.] $\vec{\alpha}$--arbitrary,\ $\beta\not=0,\;\;\vec{\gamma}=0,\;\;\vec{\phi}=0,\;\;
v=0,\;\;\vec{\xi}=0,\;\;\theta=0,\;\;\rho=0,\;\;\sigma_{ij}=0,\;\;\chi_{ij}=0,\;\;\omega_{ij}=0$\\
free parameters: $\vec{\alpha},\beta\not=0$
\item[II.] $\vec{\alpha}\not=0,\;\;\beta=0,\;\;\vec{\gamma}=0,\;\;\vec{\phi}=F\vec{\alpha},\;\;
\vec{\lambda}=L\vec{\alpha},\;\;v$--arbitrary,\ $\vec{\xi}=0,\;\;\theta=0,\;\;\rho=0,\;\;\sigma_{ij}=0,\;\;
\omega_{ij}=W(\vec{\alpha^2}\delta_{ij}-\alpha_i\alpha_j),\;\; \chi_{ij}=B(\alpha_i\alpha_j-{1\over 3}\vec{\alpha^2}\delta_{ij})+2Wv\varepsilon_{ijk}\alpha_k$\\
free parameters: $F,L,v,W\not=0,B$
\item[III.] $\vec{\alpha}=0,\;\;\beta=0,\;\;\vec{\gamma}=0,\;\;\vec{\phi}=F\vec{\mu},\;\;
\vec{\lambda}=L\vec{\mu},\;\;v=0,\;\;\vec{\xi}=0,\;\;\theta=0,\;\;\rho=0,\;\;\sigma_{ij}=0,\;\;
\omega_{ij}=W(\delta_{ij}-\mu_i\mu_j), \;\;\chi_{ij}=B(\mu_i\mu_j-{1\over 3}\delta_{ij})+C\varepsilon_{ijk}\mu_k$\\
free parameters: $F\not=0,L,B,C,W\not=0,\vec{\mu},||\vec{\mu}||=1$
\item[IV.] $\vec{\alpha}=0,\;\;\beta=0,\;\;\vec{\gamma}=0,\;\;\vec{\phi}=0,\;\;
\vec{\lambda}=L\vec{\mu},\;\;v=0,\;\;\vec{\xi}=X\vec{\mu},\;\;\theta$--arbitrary,\ $\rho=0,\;\;\sigma_{ij}=0,\;\;
\omega_{ij}=W(\delta_{ij}-\mu_i\mu_j), \;\;\chi_{ij}=B(\mu_i\mu_j-{1\over 3}\delta_{ij})+C\varepsilon_{ijk}\mu_k$\\
free parameters: $L,X,\theta,W\not=0,B,C,\vec{\mu},||\vec{\mu}||=1$
\item[V.] $\vec{\alpha}\not=0,\;\;\beta=0,\;\;\vec{\gamma}=0,\;\;\vec{\phi}=F\vec{\alpha},\;\;
\vec{\lambda}=L\vec{\alpha},\;\;v$--arbitrary,\ $\vec{\xi}=0,\;\;\theta=0,\;\;\rho=0,\;\;\sigma_{ij}=0,\;\;
\omega_{ij}=0,\;\; \chi_{ij}=B(\alpha_i\alpha_j-{1\over 3}\vec{\alpha^2}\delta_{ij})$\\
free parameters: $\vec{\alpha}\not=0,F,L,v,B$
\item[VI.] $\vec{\alpha}=0,\;\;\beta=0,\;\;\vec{\gamma}=0,\;\;\vec{\phi}=0,\;\;v=0,\;\;,\vec{\xi}=0\;\;
\theta$--arbitrary, $\;\; \rho=0,\;\;\sigma_{ij}=0,\;\;
\omega_{ij}=0,$\\
$\vec{\lambda}$\ and $\chi_{ij}$-- arbitrary except that $\varepsilon_{abc}\chi_{ab}=0$\\ 
free parameters: $\vec{\lambda},\chi_{ab},\theta$
\item[VII.] $\vec{\alpha}=0,\;\;\beta=0,\;\;\vec{\gamma}=0,\;\;\phi_i=F\varepsilon_{imn}\chi_{mn},\;\;
\vec{\lambda}=0,\;\;v=0,\;\;\vec{\xi}=0,\;\;\theta=0,\;\;\rho=0,\;\;\sigma_{ij}=0,\;\;
\omega_{ij}=0,$\\
$\chi_{ij}$-- arbitrary except that $\varepsilon_{imn}\chi_{mn}\not=0$\\ 
free parameters: $F\not=0,\chi_{mn}$
\item[VIII.] $\vec{\alpha}=0,\;\;\beta=0,\;\;\vec{\gamma}=0,\;\;\phi_i\not=F\varepsilon_{imn}\chi_{mn}\;\; \hbox{and}\;\; \vec{\phi}\not=0,\;\;
\vec{\lambda}=L\vec{\phi},\;\;v=0,\;\;\vec{\xi}=0,\;\;\theta=0,\;\;\rho=0,\;\;\sigma_{ij}=0,\;\;
\omega_{ij}=0,$\\
$\chi_{ij}$-- arbitrary\\
free parameters: $\vec{\phi}\not=0,L,\chi_{ij}$
\item[IX.] $\vec{\alpha}=0,\;\;\beta=0,\;\;\vec{\gamma}=0,\;\;\vec{\phi}\hbox{--arbitrary },\;\;
\vec{\lambda}=0,\;\;v\not=0,\;\;\vec{\xi}=0,\;\;\theta=0,\;\;\rho=0,\;\;\sigma_{ij}=0,\;\;
\omega_{ij}=0,\;\;\chi_{ij}\hbox{--arbitrary }$\\
free parameters: $\vec{\phi},v\not=0,\chi_{ij}$
\item[X.] $\vec{\alpha}=0,\;\;\beta=0,\;\;\vec{\gamma}=0,\;\;\vec{\phi}=0,\;\;
\vec{\lambda}\hbox{--arbitrary},\;\;v=0,\;\;\vec{\xi}\not=0,\;\;\theta\hbox{--arbitrary},\;\;\rho=0,\;\;\sigma_{ij}=0,\;\;
\omega_{ij}=0,\;\;\chi_{ij}\hbox{--arbitrary except that }\\ \varepsilon_{abc}\chi_{ab}\xi_c=0$\\
free parameters: $\vec{\lambda},\vec{\xi}\not=0,\theta,\chi_{ij}$
\item[XI.] $\vec{\alpha}=0,\;\;\beta=0,\;\;\vec{\gamma}=0,\;\;\vec{\phi}=0,\;\;
\vec{\lambda}=L\vec{\mu},\;\;v=0,\vec{\xi}=0,\;\;\theta\hbox{--arbitrary},\;\;\rho=-{1\over 3}S,\;\;\sigma_{ij}=S(\mu_i\mu_j-{1\over 3}\delta_{ij}),\;\;
\omega_{ij}=0,\;\;\chi_{ij}=B(\mu_i\mu_j-{1\over 3}\delta_{ij})$\\
free parameters: $S\not=0,L,B,\theta,\vec{\mu},||\vec{\mu}||=1$
\item[XII.] $\vec{\alpha}=0,\;\;\beta=0,\;\;\vec{\gamma}=0,\;\;\vec{\phi}=0,\;\;
\vec{\lambda}=0,\;\;v=0,\;\;\vec{\xi}=0,\;\;\theta\hbox{--arbitrary},\;\;\rho=-{1\over 3}S,\;\;\sigma_{ij}=S(\mu_i\mu_j-{1\over 3}\delta_{ij}),\;\;
\omega_{ij}=0,\;\;\chi_{ij}=B(\mu_i\mu_j-{1\over 3}\delta_{ij})+C\varepsilon_{ijk}\mu_k$\\
free parameters: $S\not=0,B,C\not=0,\theta,\vec{\mu},||\vec{\mu}||=1$
\item[XIII.] $\vec{\alpha}=0,\;\;\beta=0,\;\;\vec{\gamma}=0,\;\;\vec{\phi}=0,\;\;
\vec{\lambda}=0,\;\;v\not=0,\;\;\vec{\xi}=0,\;\;\theta=0,\;\;\rho=-{1\over 3}S,\;\;\sigma_{ij}=S(\mu_i\mu_j-{1\over 3}\delta_{ij}),\;\;
\omega_{ij}=0,\;\;\chi_{ij}=B(\mu_i\mu_j-{1\over 3}\delta_{ij})+C\varepsilon_{ijk}\mu_k$\\
free parameters: $S\not=0,v=0,B,C,\vec{\mu},||\vec{\mu}||=1$
\item[XIV.] $\vec{\alpha}=0,\;\;\beta=0,\;\;\vec{\gamma}=0,\;\;\vec{\phi}=0,\;\;
\vec{\lambda}=L\vec{\mu},\;\;v=0,\;\;\vec{\xi}=X\vec{\mu},\;\;\theta\hbox{--arbitrary },\;\;\rho=-{1\over 3}S,\;\;\sigma_{ij}=S(\mu_i\mu_j-{1\over 3}\delta_{ij}),\;\;
\omega_{ij}=0,\;\;\chi_{ij}=B(\mu_i\mu_j-{1\over 3}\delta_{ij})$\\
free parameters: $S\not=0,X\not=0,L,\theta,B,\vec{\mu},||\vec{\mu}||=1$
\item[XV.] $\vec{\alpha}=0,\;\;\beta=0,\;\;\vec{\gamma}=0,\;\;\vec{\phi}=F\vec{\mu},\;\;
\vec{\lambda}=0,\;\;v\not=0,\;\;\vec{\xi}=0,\;\;\theta=0,\;\;\rho=-{1\over 3}S,\;\;\sigma_{ij}=S(\mu_i\mu_j-{1\over 3}\delta_{ij}),\;\;
\omega_{ij}=0,\;\;\chi_{ij}=B(\mu_i\mu_j-{1\over 3}\delta_{ij})+C\varepsilon_{ijk}\mu_k$\\
free parameters: $S\not=0,F\not=0,v\not=0,B,C,\vec{\mu},||\vec{\mu}||=1$
\item[XVI.] $\vec{\alpha}=0,\;\;\beta=0,\;\;\vec{\gamma}=0,\;\;\vec{\phi}=F\vec{\mu},\;\;
\vec{\lambda}=0,\;\;v=0,\;\;\vec{\xi}=0,\;\;\theta=0,\;\;\rho=-{1\over 3}S,\;\;\sigma_{ij}=S(\mu_i\mu_j-{1\over 3}\delta_{ij}),\;\;
\omega_{ij}=0,\;\;\chi_{ij}=B(\mu_i\mu_j-{1\over 3}\delta_{ij})+C\varepsilon_{ijk}\mu_k$\\
free parameters: $S\not=0,F\not=0,B,C,\vec{\mu},||\vec{\mu}||=1$
\item[XVII.] $\vec{\alpha}=0,\;\;\beta=0,\;\;\vec{\gamma}=0,\;\;\vec{\phi}=F\vec{\mu},\;\;
\vec{\lambda}=L\vec{\mu},\;\;v=0,\;\;\vec{\xi}=0,\;\;\theta=0,\;\;\rho=-{1\over 3}S,\;\;\sigma_{ij}=S(\mu_i\mu_j-{1\over 3}\delta_{ij}),\;\;
\omega_{ij}=0,\;\;\chi_{ij}=B(\mu_i\mu_j-{1\over 3}\delta_{ij})$\\
free parameters: $S\not=0,F\not=0,B,L,\vec{\mu},||\vec{\mu}||=1$
\item[XVIII.] $\vec{\alpha}=0,\;\;\beta=0,\;\;\vec{\gamma}\not=0,\;\;\vec{\phi}=0,\;\;
\vec{\lambda}=L\vec{\gamma},\;\;v=0,\;\;\vec{\xi}=X\vec{\gamma},\;\;\theta=0,\;\;\rho=0,\;\;\sigma_{ij}=-\varepsilon_{ijk}\gamma_k,\;\;
\omega_{ij}=0,\;\;\chi_{ij}=0$\\
free parameters: $X,L$.
\end{description}

Let us note that all Poisson--Lie structures with $\beta=0,\;v=0$\ and $\theta=0$\ 
are coboundary and the corresponding r--matrix reads
\be
r&=&i\phi_kH\wedge P_k+i\gamma_kH\wedge K_k+i\alpha_k H\wedge J_k+\nn\\
&&+i\varepsilon_{ijk}\lambda_k P_i\wedge P_j+i(\sigma_{ij}-\rho\delta_{ij})P_i\wedge J_j\label{w20}\\
&&+i\chi_{ij} P_i\wedge K_j-i(2\omega_{ij}-\omega_{nn}\delta_{ij})J_i\wedge K_j\nn\\
&&+i\varepsilon_{ijk}\xi_kK_i\wedge K_j\nn
\ee

Now there remains only to classify the orbits under the action of residual 
automorphisms corresponding to the rotations and scaling and put our structure in the canonical
form. This is a~straightforward although very tedious task. The result can be
summarized as follows. There are 69 families of inequivalent Poisson--Lie
structures which have been grouped for convenience into eight groups. Each
group is described by an appropriate tables which are given below. They 
provide the main result of our paper.

Let us note that for all groups $\vec{n}=0$. In the last column
(labelled by \#) we indicate the number of essential parameters.
 
I. $\rho=0,\sigma_{ij}=0,\omega_{ij}=0,\chi_{ij}=0,\vec{\gamma}=0$
\vskip0.4cm
\begin{tabular}{|c|c|c|c|c|c|c|c|c|}
\hline
N&$\vec{\alpha}$&$\beta$&$\vec{\phi}$&$\vec{\lambda}$&$v$&$\vec{\xi}$&$\theta$&\#\\
\hline
1&$(0,0,\alpha)$&1&0&0&0&0&0&1\\
&$\alpha\geq0$&&&&&&&\\
2&(0,0,1)&0&(0,0,1)&(0,0,$L$)&$v$&0&0&2\\
3&(0,0,1)&0&0&(0,0,$\pm$1)&$v$&0&0&1\\
4&(0,0,1)&0&0&0&1&0&0&0\\
5&(0,0,1)&0&0&0&0&0&0&0\\
6&0&0&(0,0,1)&(0,0,$\pm$1)&0&0&0&0\\
7&0&0&(0,0,1)&0&0&0&0&0\\
8&0&0&(0,0,1)&0&1&0&0&0\\
9&0&0&0&$\lambda_1=0$&0&(0,0,1)&$\theta$&2\\
&&&&$\lambda_2^2+\lambda_3^2=1$&&&&\\
10&0&0&0&(0,0,1)&0&0&$\pm1$&0\\
11&0&0&0&(0,0,1)&0&0&0&0\\
12&0&0&0&0&0&(0,0,1)&$\theta$&1\\
13&0&0&0&0&0&0&$\pm1$&0\\
14&0&0&0&0&0&0&0&0\\
\hline
\end{tabular}
\vskip0.3cm
II. $\beta=0,\vec{\gamma}=0,\rho=0,\sigma_{ij}=0,\omega_{ij}=W(\delta_{ij}-\delta_{i3}\delta_{j3}),
\chi_{ij}=B(\delta_{i3}\delta_{j3}-{1\over 3}\delta_{ij})+C\varepsilon_{ij3}$
\vskip0.4cm
\begin{tabular}{|c|c|c|c|c|c|c|c|c|c|c|}
\hline
N&$\vec{\alpha}$&$\vec{\phi}$&$\vec{\lambda}$&$v$&$\vec{\xi}$&$\theta$&$B$&$C$&$W$&\#\\
\hline
15&(0,0,1)&(0,0,$F$)&(0,0,$L$)&$v$&0&0&$B$&$2v$&1&4\\
16&(0,0,1)&(0,0,$F$)&(0,0,$L$)&$v$&0&0&1&0&0&3\\
17&0&(0,0,$\pm1$)&(0,0,$L$)&0&0&0&$B$&$C$&1&3\\
18&0&$\vec{\phi}$&0&1&0&0&0&1&0&3\\
19&0&0&(0,0,$L$)&0&(0,0,$X$)&$\theta$&\multicolumn{2}{c|}{${\scriptstyle 2B^2+6C^2=3}$}&1&4\\
20&0&0&(0,0,$\pm1$)&0&(0,0,$X$)&$\theta$&0&0&1&2\\
21&0&0&0&0&(0,0,$X$)&$\theta$&0&0&1&2\\
\hline
\end{tabular}
\vskip0.5cm
III. $\vec{\alpha}=0,\beta=0,\vec{\gamma}=0,\rho=-{1\over 3},
\sigma_{ij}=(\delta_{i3}\delta_{j3}-{1\over 3}\delta_{ij}),\omega_{ij}=0,
\chi_{ij}=B(\delta_{i3}\delta_{j3}-{1\over 3}\delta_{ij})+C\varepsilon_{ij3}$
\vskip0.4cm
\begin{tabular}{|c|c|c|c|c|c|c|c|c|}
\hline
N&$\vec{\phi}$&$\vec{\lambda}$&$v$&$\vec{\xi}$&$\theta$&$B$&$C$&\#\\
\hline
22&(0,0,1)&(0,0,$L$)&0&0&0&$B$&0&2\\
23&(0,0,1)&0&$v\not=0$&0&0&$B$&0&2\\
24&(0,0,1)&0&$v$&0&0&$B$&$C\not=0$&3\\
25&0&(0,0,$L$)&0&(0,0,$\pm1$)&$\theta$&$B$&0&3\\
26&0&(0,0,$L$)&0&0&$\pm1$&$B$&0&2\\
27&0&(0,0,$L$)&0&0&0&1&0&1\\
28&0&(0,0,$L$)&0&0&0&0&0&1\\
29&0&0&1&0&0&$B$&$C$&2\\
30&0&0&0&0&$\theta$&\multicolumn{2}{l|}{$\scriptstyle 2B^2+6C^2=3$}&2\\
&&&&&&\multicolumn{2}{c|}{$\scriptstyle C\not=0$}&\\
\hline
\end{tabular}

\vskip0.5cm
IV. $\vec{\alpha}=0,\beta=0,\vec{\phi}=0,v=0,,\theta=0,\rho=0,
\omega_{ij}=0,
\chi_{ij}=0,\sigma_{ij}=-\varepsilon_{ij3}$
\vskip0.4cm
\begin{tabular}{|c|c|c|c|c|}
\hline
N&$\vec{\gamma}$&$\vec{\lambda}$&$\vec{\xi}$&\#\\
\hline
31&(0,0,1)&(0,0,$L$)&(0,0,$\pm1$)&1\\
32&(0,0,1)&(0,0,$L$)&0&1\\
\hline
\end{tabular}
\vskip0.5cm
V.$\vec{\alpha}=0,\beta=0,\vec{\gamma}=0,v=0,\rho=0,\sigma_{ij}=0,
\omega_{ij}=0,tr\chi=0,\sum\limits^3_{ij}\chi_{ij}^2=1,\chi_{ij}=diag(\chi_{11},\chi_{22},\chi_{33})$

$\chi_{11}\not=\chi_{22}\not=\chi_{33}$
\vskip0.4cm
\begin{tabular}{|c|c|c|c|c|c|}
\hline
N&$\vec{\phi}$&$\vec{\lambda}$&$\vec{\xi}$&$\theta$&\#\\
\hline
33&0&$||\vec{\lambda}||=1$&0&$\theta$&4\\
34&0&$\vec{\lambda}$&$||\vec{\xi}||=1$&$\theta$&7\\
35&$||\vec{\phi}||=1$&$\vec{\lambda}=L\vec{\phi}$&0&0&4\\
\hline
\end{tabular}
\vskip0.4cm
(Vb)\hskip0.4cm$\chi_{11}=\chi_{22}\not=\chi_{33}$
\vskip0.3cm
\begin{tabular}{|c|c|c|c|c|c|}
\hline
N&$\vec{\phi}$&$\vec{\lambda}$&$\vec{\xi}$&$\theta$&\#\\
\hline
36&0&$\lambda_1=0$&0&$\theta$&2\\
&&$\lambda_2^2+\lambda_3^2=1$&&&\\
37&0&$\vec{\lambda}$&$\xi_1=0$&$\theta$&5\\
&&&$\xi_2^2+\xi_3^2=1$&&\\
38&$\phi_1=0$&$\vec{\lambda}=L\vec{\phi}$&0&0&2\\
&$\phi_2^2+\phi_3^2=1$&&&&\\
\hline
\end{tabular}
\vskip0.4cm

VI.$\vec{\alpha}=0,\beta=0,\vec{\gamma}=0,v=0,\rho=0,\sigma_{ij}=0,
\omega_{ij}=0,tr\chi=0,\sum\limits^3_{ij}\chi_{ij}^2=1,\chi_{32}=-\chi_{23}\not=0,
\chi_{12}=\chi_{21},\chi_{13}=\chi_{32}$

(VIa)\hskip0.4cm$\chi_{22}\not=\chi_{33}$
\vskip0.3cm
\begin{tabular}{|c|c|c|c|c|c|}
\hline
N&$\vec{\phi}$&$\vec{\lambda}$&$\vec{\xi}$&$\theta$&\#\\
\hline
39&0&$\lambda_1=0$&0&$\theta$&6\\
&&$\lambda_2^2+\lambda_3^2=1$&&&\\
40&0&$\vec{\lambda}$&$\xi_1=0$&$\theta$&9\\
&&&$\xi_2^2+\xi_3^2=1$&&\\
41&$||\vec{\phi}||=1$&$\vec{\lambda}=L\vec{\phi}$&0&0&7\\
&$\phi_2^2+\phi_3^2=1$&&&&\\
\hline
\end{tabular}
\vskip0.5cm
(VIb)\hskip0.4cm$\chi_{22}=\chi_{33}$
\vskip0.3cm
\begin{tabular}{|c|c|c|c|c|c|}
\hline
N&$\vec{\phi}$&$\vec{\lambda}$&$\vec{\xi}$&$\theta$&\#\\
\hline
42&0&$\lambda_1=\lambda_2=0$&0&$\theta$&4\\
&&$\lambda_2=\pm1$&&&\\
43&0&$\vec{\lambda}$&$\xi_1=\xi_2=0$&$\theta$&7\\
&&&$\xi_3=\pm1$&&\\
44&($F$,0,1)&$\vec{\lambda}=L\vec{\phi}$&0&0&5\\
\hline
\end{tabular}
\vskip0.5cm
VII.$\vec{\alpha}=0,\beta=0,\vec{\gamma}=0,\vec{\lambda}=0,\vec{\xi}=0,
\rho=0,\sigma_{ij}=0,
\omega_{ij}=0,tr\chi=0,\sum\limits^3_{ij}\chi_{ij}^2=1,$
\vskip0.4cm
(VIIa)\hskip0.4cm$\chi_{11}\not=\chi_{22}\not=\chi_{33},\chi_{13}=-\chi_{31},\chi_{23}=-\chi_{32},\chi_{12}=-\chi_{21}$
\vskip0.3cm
\begin{tabular}{|c|c|c|c|c|}
\hline
N&$\vec{\phi}$&$v$&$\theta$&\#\\
\hline
45&0&0&$\pm1,0$&4\\
46&$\vec{\phi}$&1&0&7\\
\hline
\end{tabular}
\vskip0.4cm
(VIIb)\hskip0.4cm$\chi_{11}=\chi_{22}\not=\chi_{33},\chi_{13}=-\chi_{31}\not=0,\chi_{23}=\chi_{32}=0,\chi_{12}=-\chi_{21}$
\vskip0.3cm
\begin{tabular}{|c|c|c|c|c|}
\hline
N&$\vec{\phi}$&$v$&$\theta$&\#\\
\hline
47&0&0&$\pm1,0$&2\\
48&$\vec{\phi}$&1&0&5\\
\hline
\end{tabular}
\vskip0.4cm
(VIIc)\hskip0.4cm$\chi_{11}=\chi_{22},\chi_{13}=\chi_{31}=\chi_{23}=\chi_{32}=0,\chi_{12}=-\chi_{21}$
\vskip0.3cm
\begin{tabular}{|c|c|c|c|c|}
\hline
N&$\vec{\phi}$&$v$&$\theta$&\#\\
\hline
49&0&0&$\pm1,0$&1\\
50&$(0,\phi_2,\phi_3)$&1&0&3\\
\hline
\end{tabular}
\vskip0.5cm
VIII.$\vec{\alpha}=0,\beta=0,\vec{\gamma}=0,\vec{\phi}=(0,0,1),
\vec{\lambda}=0,v=0,\vec{\xi}=0,\theta=0,
\rho=0,\sigma_{ij}=0,
\omega_{ij}=0,tr\chi=0,\sum\limits^3_{ij}\chi_{ij}^2=1,\chi_{12}=-\chi_{21}={1\over F}$
\vskip0.3cm
(VIIIa)\hskip0.4cm$\chi_{11}\not=\chi_{22},\chi_{13}=\chi_{31},\chi_{23}=\chi_{32}$

(VIIIb)\hskip0.4cm$\chi_{11}=\chi_{22},\chi_{13}=\chi_{31}=0,\chi_{23}=\chi_{32}$
\\

Now, inserting the appropriate values of the 
parameters (listed above) to eq.\mref{wb1} 
and using eq.\mref{w19} one can easily calculate the fundamental Poisson brackets
for all nonequivalent structures.

\section{Summary}

We have obtained all Lie--Poisson structures
on the four-dimensional Galilei group and classified
them up to the equivalence implied by group
automorphisms. The resulting set of structures appears
to be quite rich; in particular, it includes many
non--coboundary structures, to be contrasted with the
Poincare group case~\cite{b5}. In spite of that, part of them
can surely be obtained from those on the Poincare group
by a~contraction procedure.

The next step to be done is to quantize the 
Lie-Poisson structures. In general, the consistent 
quantization is not an easy task. However, the preliminary
study already done by us shows, that most of the 
cases described here are quantization friendly.

\section{Appendix A: The cocycle condition for subgroups}
\renewcommand{\theequation}{A.\arabic{equation}}
\zero
\indent Let us specify the equations \mref{w17} for the subgroups of
rotations, boosts, space and time translations. They read, respectively:
\be
\Psi_i(RR')&=&\Psi_i(R)+R_{il}\Psi_l(R')\nn\\
\Phi_i(RR')&=&\Phi_i(R)+R_{il}\Phi_l(R')\nn\\
\Gamma_i(RR')&=&\Gamma_l(R)+R_{il}\Gamma_l(R')\nn\\
\Lambda_i(RR')&=&\Lambda_i(R)+R_{il}\Lambda_l(R')\nn\\
\Upsilon_{ij}(RR')&=&\Upsilon_{ij}(R)+R_{im}R_{jl}\Upsilon_{ml}(R')\label{wa1a}\\
\Sigma_{ij}(RR')&=&\Sigma_{ij}(R)+R_{im}R_{jl}\Sigma_{ml}(R')\nn\\
\Xi_{i}(RR')&=&\Xi_{i}(R)+R_{im}\Xi_{m}(R')\nn\\
\Omega_{ij}(RR')&=&\Omega_{ij}(R)+R_{im}R_{jl}\Omega_{ml}(R')\nn\\
\Pi_i(RR')&=&\Pi_i(R)+R_{il}\Pi_l(R')\nn\\
&&\nn\\
\Psi_i(\vec{v}+\vec{v'})&=&\Psi_i(\vec{v})+\Psi_i(\vec{v'})\nn\\
\Phi_i(\vec{v}+\vec{v'})&=&\Phi_i(\vec{v})+\Phi_i(\vec{v'})\nn\\
\Gamma_i(\vec{v}+\vec{v'})&=&\Gamma_i(\vec{v})+\Gamma_i(\vec{v'})-\varepsilon_{ink}v_n\Psi_k(\vec{v'})\nn\\
\Lambda_i(\vec{v}+\vec{v'})&=&\Lambda_i(\vec{v})+\Lambda_i(\vec{v'})-\fr\varepsilon_{imn}v_m\Phi_n(\vec{v'})\label{wa1b}\\
\Upsilon_{ij}(\vec{v}+\vec{v'})&=
&\Upsilon_{ij}(\vec{v})+\varepsilon_{jnk}v_n\Psi_k(\vec{v'})v_i-
\Gamma_j(\vec{v'})v_i-\varepsilon_{jnk}v_n\Sigma_{ik}(\vec{v'})\nn\\
\Sigma_{ij}(\vec{v}+\vec{v'})&=&\Sigma_{ij}(\vec{v})+\Sigma_{ij}(\vec{v'})-v_i\Psi_j(\vec{v'})\nn\\
\Xi_i(\vec{v}+\vec{v'})&=&\Xi_i(\vec{v})+\Xi_i(\vec{v'})-\fr(v_i\Omega_{mm}(\vec{v'})-\nn\\
&&-v_p\Omega_{pi}(\vec{v'}))+\Pi_m(\vec{v'})v_mv_i\nn\\
\Omega_{ij}(\vec{v}+\vec{v'})&=&\Omega_{ij}(\vec{v})+\Omega_{ij}(\vec{v'})+2\Pi_i(\vec{v'})v_j-2\Pi_m(\vec{v'})v_m\delta_{ij}\nn\\
\Pi_i(\vec{v}+\vec{v'})&=&\Pi_i(\vec{v})+\Pi_i(\vec{v'})\nn\\
&&\nn\\
\Psi_i(\vec{a}+\vec{a'})&=&\Psi_i(\vec{a})+\Psi_i(\vec{a'})\nn\\
\Phi_i(\vec{a}+\vec{a'})&=&\Phi_i(\vec{a})+\Phi_i(\vec{a'})-\varepsilon_{ink}a_n\Psi_k(\vec{a'})\nn\\
\Gamma_i(\vec{a}+\vec{a'})&=&\Gamma_i(\vec{a})+\Gamma_i(\vec{a'})\nn\\
\Lambda_i(\vec{a}+\vec{a'})&=&\Lambda_i(\vec{a})+\Lambda_i(\vec{a'})+\Pi_m(\vec{a'})a_ma_i-\fr(a_i\delta_{lm}-a_l\delta_{im})\Sigma_{lm}(\vec{a'})\nn\\
\Upsilon_{ij}(\vec{a}+\vec{a'})&=&\Upsilon_{ij}(\vec{a})+\Upsilon_{ij}(\vec{a'})+\varepsilon_{ink}a_n\Omega_{jk}(\vec{a'})\label{wa1c}\\
\Sigma_{ij}(\vec{a}+\vec{a'})&=&\Sigma_{ij}(\vec{a})+\Sigma_{ij}(\vec{a'})+2(\Pi_i(\vec{a'})a_j-\Pi_l(\vec{a'})a_l\delta_{ij})\nn\\
\Xi_{i}(\vec{a}+\vec{a'})&=&\Xi_i(a)+\Xi_i(a)\nn\\
\Omega_{ij}(\vec{a}+\vec{a'})&=&\Omega_{ij}(\vec{a})+\Omega_{ij}(\vec{a'})\nn\\
\Pi_i(\vec{a}+\vec{a'})&=&\Pi_i(\vec{a})+\Pi_i(\vec{a'})\nn\\
&&\nn\\
\Psi_i({t}+{t'})&=&\Psi_i({t})+\Psi_i({t'})\nn\\
\Phi_i({t}+{t'})&=&\Phi_i({t})+\Phi_i({t'})-t\Gamma_i(t')\nn\\
\Gamma_i({t}+{t'})&=&\Gamma_i({t})+\Gamma_i({t'})\label{wa1d}\\
\Lambda_i({t}+{t'})&=&\Lambda_i({t})+\Lambda_i({t'})-\fr t\varepsilon_{imn}\Upsilon_{mn}(t')+t^2\Xi_i(t')\nn\\
\Upsilon_{ij}(t+t')&=&\Upsilon_{ij}(t)+\Upsilon_{ij}(t')-2t\varepsilon_{ijl}\Xi_l(t')\nn\\
\Sigma_{ij}(t+t')&=&\Sigma_{ij}(t)+\Sigma_{ij}(t')-t\Omega_{ij}(t')\nn\\
\Xi_i(t+t')&=&\Xi_i(t)+\Xi_i(t')\nn\\
\Omega_{ij}(t+t')&=&\Omega_{ij}(t)+\Omega_{ij}(t')\nn\\
\Pi_i(t+t')&=&\Pi_i(t)+\Pi_i(t')\nn
\ee
Note that all eqs\mref{wa1a} have the same structure: 
\be
T_{i_1\ldots i_k}(RR')&=&T_{i_1\ldots i_k}(R)+
R_{i_1j_1}\ldots R_{i_kj_k}T_{j_1\ldots j_k}(R')\label{wa2}
\ee
They can be solved by integrating over $R'$\ with respect to the
Haar measure on SO(3):
\be
T_{i_1\ldots i_k}(R)&=&(R_{i_1j_1}\ldots R_{i_kj_k}-\delta_{i_1 j_1}\ldots\delta_{i_kj_k})c_{j_1\ldots j_k}\label{wa3a}
\ee
where $c_{j_1\ldots j_k}$ are constants.
This result agrees with general theorem that all cocycles on semisimple groups
are coboundaries. On the other hand it follows immediately from \hbox{eqs.(\ref{wa1b}--\ref{wa1d})}
that all functions entering there are polynomials in the relevant parameters. 
This allows us to write out explicitly the general solutions.
\be
\Psi_i(\vec{v})&=&0\nn\\
\Phi_i(\vec{v})&=&0\nn\\
\Gamma_i(\vec{v})&=&a_{ij}v_j\nn\\
\Lambda_i(\vec{v})&=&b_{ij}v_j\label{wa3b}\\
\Upsilon_{ij}(\vec{v})&=&c_{ijk}v_k-\fr(\delta_{ik}a_{jl}+\delta_{ji}a_{lk}-a_{lj}\delta_{ki}+\fr\varepsilon_{jkl}\varepsilon_{inm}a_{nm})v_kv_l\nn\\
\Sigma_{ij}(\vec{v})&=&(\fr\delta_{jk}\varepsilon_{inm}a_{nm}-\varepsilon_{ijn}a_{kn})v_k\nn\\
\Xi_i(\vec{v})&=&d_{ik}v_k+({1\over 4}e_{jik}-{1\over 8}(e_j\delta_{ik}+e_k\delta_{ij}))v_jv_k\nn\\
\Omega_{ij}(\vec{v})&=&(e_{ijk}+\fr(e_k\delta_{ji}-e_i\delta_{jk}))v_k\nn\\
\Pi_i(\vec{v})&=&0,\nn
\ee
where $e_{ijk}=e_{kji},\;e_{iik}=0$;
\be
\Psi_i(\vec{a})&=&0\nn\\
\Phi_i(\vec{a})&=&f_{ij}a_j\nn\\
\Gamma_i(\vec{a})&=&g_{ij}a_j\nn\\
\Lambda_i(\vec{a})&=&h_{ij}a_j+{1\over 4}(h_{jik}-\fr(h_j\delta_{ik}+h_k\delta_{ij}))a_ja_k\nn\\
\Upsilon_{ij}(\vec{a})&=&k_{ijk}a_k\label{wa3c}\\
\Sigma_{ij}(\vec{a})&=&(h_{ijk}+\fr(h_k\delta_{ij}-h_i\delta_{jk}))a_k\nn\\
\Xi_i(\vec{a})&=&l_{ij}a_j\nn\\
\Omega_{ij}(\vec{a})&=&0\nn\\
\Pi_i(\vec{a})&=&0,\nn
\ee
where $h_{ijk}=h_{kji},\;h_{iik}=0$;
\be
\Psi(t)&=&p_it\nn\\
\Phi_i(t)&=&r_it-\fr s_it^2\nn\\
\Gamma_i(t)&=&s_it\nn\\
\Lambda_i(t)&=&u_it-{1\over 4}\varepsilon_{imn}x_{mn}t^2+{1\over 3}w_it^3\nn\\
\Upsilon_{ij}(t)&=&x_{ij}t-\varepsilon_{ijk}w_kt^2\label{wa3d}\\
\Sigma_{ij}(t)&=&y_{ij}t-\fr z_{ij}t^2\nn\\
\Xi_i(t)&=&w_it\nn\\
\Omega_{ij}(t)&=&z_{ij}t\nn\\
\Pi_i(t)&=&m_it.\nn
\ee
\section{Appendix B: The general solution to the cocycle condition}
\renewcommand{\theequation}{B.\arabic{equation}}
\zero
\indent According to the procedure outlined in sec.\ref{s2}, we use
the cocycle condition(i) together with the decomposition \mref{w18} and
the expressions written out in Appendix~A to produce the Ansatz
for \hbox{$\eta(g)$.} Inserting it back into \mref{w17} we find the
general solution for $\eta$\ of the form.
\be
\Psi_i(g)&=&(R_{ij}-\delta_{ij})\alpha_j\nn\\
\Phi_i(g)&=&(R_{ij}-\delta_{ij})\phi_j+\beta(a_i-v_it)-\gamma_jR_{ij}t+\varepsilon_{ijk}
\alpha_lR_{jl}(a_k-v_kt)\nn\\
\Gamma_i(g)&=&(R_{ij}-\delta_{ij})\gamma_j+\beta v_i+\varepsilon_{ijk}\alpha_lR_{jl}v_k\nn\\
\Lambda_i(g)&=&(R_{ij}-\delta_{ij})\lambda_j+(\rho-\fr\sigma_{nn})(a_i-v_i t)+\fr\beta\varepsilon_{ijk}a_jv_k\label{wb1}\\
&&+\fr\varepsilon_{ijk}\phi_lR_{jl}v_k-\fr\alpha_jR_{ij}(a_kv_k-\vec{v}^2t)+\nn\\
&&+\fr\alpha_jR_{kj}v_k(a_i-v_it)+\xi_jR_{ij}t^2-\fr\varepsilon_{ijk}\gamma_lR_{jl}v_kt+\nn\\
&&-\fr\varepsilon_{jkl}\chi_{kl}R_{ij}t-\omega_{jl}R_{ij}R_{kl}a_kt+\nu v_i+\nn\\
&&+\fr\sigma_{ij}R_{ij}R_{kl}(a_k-v_kt)+n_kR_{mk}(a_ia_m+2v_mv_it)-\nn\\
&&-n_l(a_kR_{kl}v_i+a_iR_{kl}v_k)t\nn\\
\Upsilon_{ij}(g)&=&(R_{ik}R_{jl}-\delta_{ik}\delta_{jl})\chi_{kl}+\delta_{ij}\theta t-
\fr\delta_{ij}\beta\vec{v}^2+\nn\\
&&+\varepsilon_{jkl}\alpha_nR_{ln}v_kv_i-\gamma_kR_{jk}v_i-2\varepsilon_{ijk}\xi_lR_{kl}t\nn\\
&&-\varepsilon_{jkl}\sigma_{ns}R_{in}R_{ls}v_k+\rho\varepsilon_{ijk}v_k+\omega_{nn}\varepsilon_{ijk}a_k\nn\\
&&+2\omega_{ns}(\varepsilon_{ikl}R_{js}R_{ln}a_k-\varepsilon_{ijl}R_{ks}R_{ln}v_kt)\nn\\
&&+2n_l(R_{jl}v_k-R_{ml}v_m\delta_{ik})\varepsilon_{kin}a_n+\nn\\
&&-2n_kv_sv_mR_{mk}\varepsilon_{sij}t+\omega_{lk}R_{nk}R_{il}v_nt^2\nn\\
\Sigma_{ij}(g)&=&(R_{ik}R_{jl}-\delta_{ik}\delta_{jl})\sigma_{kl}-\beta\varepsilon_{ijk}v_k-\alpha_kR_{jk}v_i\nn\\
&&-2\omega_{lk}R_{ik}R_{jl}t+\omega_{nn}\delta_{ij}t+\nn\\
&&+2n_k(R_{ik}a_j-R_{mk}a_m\delta_{ij}-R_{ik}v_jt+R_{mk}v_m\delta_{ij}t)\nn\\
\Xi_i(g)&=&(R_{ij}-\delta_{ij})\xi_j+\omega_{jl}R_{ij}R_{kl}v_k+n_kR_{mk}v_mv_i\nn\\
\Omega_{ij}(g)&=&2(R_{ik}R_{jl}-\delta_{ik}\delta_{jl})\omega_{lk}+2n_k(R_{ik}v_j-R_{mk}v_m\delta_{ij})\nn\\
\Pi_i(g)&=&(R_{ij}-\delta_{ij})n_j\nn
\ee
\section{Appendix C: Jacobi identities}
\renewcommand{\theequation}{C.\arabic{equation}}
\zero
Using the general form of $\eta$\ described in Appendix~B we find from eq.\mref{w13}:
\be
\delta(H)&=&\gamma_iH\wedge P_i+\fr(\chi_{kj}-\chi_{jk})P_j\wedge P_k+(2\varepsilon_{ijk}\xi_k-\theta\delta_{ij})P_i\wedge K_j\nn\\
&&(2\omega_{ji}-\omega_{nn}\delta_{ji})P_i\wedge J_j\nn\\
\delta(P_s)&=&(\beta\delta_{is}+\varepsilon_{iks}\alpha_k)H\wedge P_i+\varepsilon_{ijk}(\rho\delta_{is}-\fr\sigma_{nn}\delta_{is}+\nn\\
&&+\fr\sigma_{si})P_j\wedge P_k+(2\varepsilon_{lis}\omega_{lj}-\varepsilon_{jis}\omega_{nn})P_i\wedge K_j\nn\\
&&+2(n_i\delta_{sj}-n_s\delta_{ij})P_i\wedge J_j\label{wc1}\\
\delta(K_s)&=&(\beta\delta_{is}+\varepsilon_{iks}\alpha_k)H\wedge P_i+(v\varepsilon_{sjk}-\fr(\phi_k\delta_{js}-\phi_j\delta_{ks}))P_j\wedge P_k\nn\\
&&+(\rho\varepsilon_{ijs}-\varepsilon_{kjs}\sigma_{ik}-\gamma_j\delta_{is})P_i\wedge K_j+\nn\\
&&-(\beta\varepsilon_{ijs}+\alpha_j\delta_{is})P_i\wedge J_j+\varepsilon_{ijk}\omega_{is}K_j\wedge K_k+\nn\\
&&+2(n_j\delta_{sk}-n_s\delta_{jk})K_j\wedge J_k\nn\\
\delta(J_s)&=&\varepsilon_{sij}\alpha_jH\wedge J_i+\varepsilon_{sij}\phi_jH\wedge P_i+\varepsilon_{ijs}\gamma_j H\wedge K_i\nn\\
&&+(\lambda_j\delta_{ks}-\lambda_k\delta_{js})P_j\wedge P_k+(\varepsilon_{sik}\lambda_{kj}+\varepsilon_{sjk}\chi_{ik})P_i\wedge K_j\nn\\
&&+(\varepsilon_{sik}\sigma_{kj}+\varepsilon_{sjk}\sigma_{ik})P_i\wedge J_j+(\xi_i\delta_{ks}-\xi_k\delta_{js})K_j\wedge K_k\nn\\
&&+2(\varepsilon_{sik}\omega_{jk}+\varepsilon_{sjk}\omega_{ki})K_i\wedge J_j+2n_kJ_k\wedge J_s\nn
\ee
Let $\tilde{X_i}$\ denote a basis in ${\cal G}^*$\ defined by $<\tilde{X_i},X_j>=\delta_{ij}$.
Then $\delta$\ imposes the following commutator structure in ${\cal G}^*$.
\be
{}[\tilde{H},\tilde{J_k}]&=&\varepsilon_{ikl}\alpha_l\tilde{J_i}\nn\\
{}[\tilde{H},\tilde{P_k}]&=&\gamma_k\tilde{H}+(\beta\delta_{ik}+\varepsilon_{ikl}\alpha_l)\tilde{P_i}+\varepsilon_{ikl}\phi_l\tilde{J_i}\nn\\
{}[\tilde{H},\tilde{K_k}]&=&(\beta\delta_{ik}+\varepsilon_{ikl}\alpha_l)\tilde{K_i}+\varepsilon_{ikl}\gamma_l\tilde{J_i}\nn\\
{}[\tilde{K_k},\tilde{J_l}]&=&2(n_k\delta_{li}-n_i\delta_{kl})\tilde{K_i}+2(\varepsilon_{ikn}\omega_{ln}+\varepsilon_{iln}\omega_{nk})\tilde{J_i}\nn\\
{}[\tilde{P_l},\tilde{P_m}]&=&(\chi_{lm}-\chi_{ml})\tilde{H}+2\varepsilon_{klm}[\rho\delta_{ki}+\fr(\sigma_{ik}-\sigma_{nn}\delta_{ik})]\tilde{P_i}\nn\\
&&+2[v\varepsilon_{ilm}+\fr(\phi_l\delta_{im}-\phi_m\delta_{li})]\tilde{K_i}+\label{wc2}\\
&&+2(\lambda_l\delta_{im}-\lambda_m\delta_{ik})\tilde{J_i}\nn\\
{}[\tilde{P_k},\tilde{K_l}]&=&(2\xi_n\varepsilon_{nkl}-\theta\delta_{kl})\tilde{H}+(2\varepsilon_{nki}\omega_{nl}-\omega_{nn}\varepsilon_{lki})\tilde{P_i}\nn\\
&&+(\rho\varepsilon_{kli}-\varepsilon_{lin}\sigma_{kn}-\delta_{ki}\gamma_l)\tilde{K_i}+\nn\\
&&+(\varepsilon_{ikn}\chi_{nl}+\varepsilon_{iln}\chi_{kn})\tilde{J_i}\nn\\
{}[\tilde{P_k},\tilde{J_l}]&=&(2\omega_{ik}-\omega_{nn}\delta_{lk})\tilde{H}+2(n_k\delta_{li}-n_i\delta_{kl})\tilde{P_i}\nn\\
&&-(\beta\varepsilon_{kli}+\alpha_l\delta_{ki})\tilde{K_i}+(\varepsilon_{ikn}\sigma_{nl}+\varepsilon_{iln}\sigma_{kn})\tilde{J_i}\nn\\
{}[\tilde{K_m},\tilde{K_n}]&=&2\varepsilon_{kmn}\omega_{ki}\tilde{K_i}+2(\xi_m\delta_{ni}-\xi_n\delta_{mi})\tilde{J_i} \nn\\
{}[\tilde{J_k},\tilde{J_l}]&=&2(n_k\delta_{li}-n_l\delta_{ki})\tilde{J_i}\nn
\ee
Now we have to solve the Jacobi identities for the structure described in \mref{wc2}.
This is still a~complicated task so we apply a~mixed procedure consisting in 
solving part of Jacobi identities directly and applying the group of automorphisms
to simplify the remaining ones.

Therefore, we first give the action of automorphism group on parameters
of $\delta$. Let us start with  inner automorphisms. They  read
\be
\vec{\alpha_k}=\vec{\alpha_i}&&\gamma'_k=\gamma_k-\beta v_k+\varepsilon_{kil}v_i\alpha_l\nn\\
\beta'=\beta&&\lambda_k'=\lambda_k-vv_k-\fr\varepsilon_{kil}\phi_iv_\alpha\nn\\
\vec{\phi'}=\vec{\phi}&&\xi'_k=\xi_k-\omega_{kn}v_n-v_k(v_nn_n)\label{wc3a}\\
\theta'=\theta&&\rho'=\rho-{1\over 3}\alpha_nv_n\nn\\
v'=v&&\omega'_{ij}=\omega_{ij}+v_in_j+(v_nn_n)\delta_{ij}\nn\\
\vec{n'}=\vec{n}&&\sigma'_{ia}=\sigma_{ia}+\beta\varepsilon_{ial}v_l+\alpha_av_i-{1\over 3}\alpha_nv_n\delta_{ai}\nn\\
&&\chi'_{ab}=\chi_{ab}+v_a\gamma_b+\varepsilon_{bnm}v_n\alpha_mv_a-\fr\varepsilon_{abk}\sigma_{nk}v_n\nn\\
&&-\rho\varepsilon_{abn}v_n-\fr\varepsilon_{alm}\sigma_{bl}v_m-\fr\varepsilon_{blm}\sigma_{al}v_m\nn\\
&&-{1\over 3}(v_n\gamma_n+\varepsilon_{lnm}\sigma_{nl}v_m)\delta_{ab}\nn
\ee
for boosts,
\be
\vec{n'}=\vec{n},&\vec{\alpha'}={\alpha},&\beta'=\beta\nn\\
v'=v,&\vec{\xi'}=\vec{\xi},&\theta'=\theta\nn
\ee
\be
\phi_i'&=&\phi_i-\beta a_i-\varepsilon_{ink}\alpha_na_k\nn\\
\gamma'_k&=&\gamma_k-2\varepsilon_{knm}a_mn_n\label{wc3b}\\
\lambda'_k&=&\lambda_k-\fr\sigma_{nk}a_n-\rho a_k-\vec{a^2}n_k-(n_ia_i)n_k\nn\\
\rho'&=&\rho-{4\over 3}a_kn_k\;\;\sigma'_{ab}=\sigma_{ab}+{2\over 3}(a_kn_k)\delta_{ab}-2n_aa_b\nn\\
\chi'_{ab}&=&\chi_{ab}-\varepsilon_{abn}\omega_{nk}a_k+\varepsilon_{amk}\omega_{mb}a_k+\nn\\
&&+\varepsilon_{bmk}\omega_{ma}a_k+{2\over 3}\varepsilon_{mnk}\omega_{mn}a_k\delta_{ab}\nn\\
\omega'_{ab}&=&\omega_{ab}\nn
\ee
for space translations and 
\be
\vec{n'}=\vec{n},&\vec{\alpha}=\alpha,&\beta'=\beta,\;\;\;\;\vec{\gamma'}=\vec{\gamma}\nn\\
v'=v,&\vec{\xi'}=\vec{\xi},&\theta'=\theta,\;\;\;\;\omega'_{ij}=\omega_{ij}\nn\\
\vec{\phi'}=\vec{\phi}+t\vec{\gamma}&,&\lambda'_i=\lambda_i-\fr t\varepsilon_{imn}\chi_{nm}+t^2\xi_i\nn\\
\rho'=\rho+{1\over 3}t\omega_{nn}&,&\sigma'_{ab}=\sigma_{ab}+2t\omega_{ba}-{2\over 3}t\omega_{nn}\delta_{ab}\nn\\
\chi'_{ab}=\chi_{ab}+2t\varepsilon_{abk}\xi_k\nn
\ee
for time translations, respectively.

Under the rotations the parameters transform as tensors of appropriate rank.

Besides, there are two outer automorphisms, which correspond to rescaling
of space and time unit, \hbox{($\vec{a}\to a\vec{a},\;t\to bt$)}. They read:
\be
n'=\vec{n},&\alpha'={1\over b}\vec{\alpha},&\beta'={1\over b}\beta\nn\\
\vec{\gamma}={1\over a}\gamma,&\vec{\phi}={1\over ab}\vec{\phi},&\vec{\lambda}={1\over a^2}\vec{\lambda}\nn\\
v'={1\over ab}v,&\vec{\xi}={b^2\over a^2}\vec{\xi},&\theta'={b^2\over a^2}\theta\label{wc4}\\
\rho'={1\over a}\rho,&\sigma'_{ab}={1\over a}\sigma_{ab},&\chi'_{ab}={b\over a^2}\chi_{ab}\nn\\
\omega'_{ab}={b\over a}\omega_{ab}\nn
\ee
We shall not enter into all details. Let us rather give a sketch of the procedure.

By solving Jacobi identities for the subalgebra generated by $\tilde{H},\tilde{J_l},\tilde{K_n}$\ 
we find the following six families of constraints on parameters $\vec{n},\beta,\vec{\alpha},
\omega_{ij},\vec{\gamma},\hbox{\ and\ }\vec{\xi}$.\\
a)$\vec{n}=0,\vec{\alpha}\hbox{--arbitrary\ },\beta\not=0,\vec{\gamma}\hbox{--arbitrary,\ }\vec{\xi}=0,
\omega_{ij}=0$\\
b)$\vec{n}=0,\vec{\alpha}\not=0,\beta=0,\vec{\gamma}\hbox{--arbitrary,\ }\vec{\xi}=\xi\vec{\alpha}+W(\vec{\alpha}\times\vec{\gamma}),
\omega_{ij}=W(\vec{\alpha^2}\delta_{ij}-\alpha_i\alpha_j), W\not=0$\\
c)$\vec{n}=0,\vec{\alpha}=0,\beta=0,\vec{\xi}\hbox{--arbitrary,\ }\omega_{ij}=W(\delta_{ij}-\mu_i\mu_j)+V\varepsilon_{ijk}\mu_k,
||\vec{\mu}||=1$\\
$\vec{\gamma}=\left\{
\begin{array}{c}
\gamma\vec{\mu}\;\;\;\;\; \hbox{if\ } |W|+|V|\not=0\\
\hbox{arbitrary\ }\;\hbox{if\ } W=V=0
\end{array}\right.$\\
d)$\vec{n}=0,\vec{\alpha}\not=0,\beta=0,\omega_{ij}=W\delta_{ij},\vec{\gamma}={1\over W}(\vec{\alpha}\times\vec{\xi}),
W\not=0,\vec{\xi}\hbox{--arbitrary}$\\
e)$\vec{n}=0,\vec{\alpha}\not=0,\beta=0,\omega_{ij}=0,\vec{\xi}||\vec{\alpha},\vec{\gamma}\bot\vec{\alpha}$\\
f)$\vec{n}=0,\vec{\alpha}=0,\beta=0,\omega_{ij}=W\delta_{ij},\vec{\gamma}=0,\vec{\xi}\hbox{--arbitrary}$\\

Now we have to solve the remaining Jacobi identities. First of all let us
note that there are ambiguities in determining the matrices $\sigma_{ij}$\ and
$\chi_{ij}$. Namely, both can be redefined by adding the arbitrary multiplies
of unit matrix. In order to remove this ambiguity we put 
\hbox{$tr\sigma=tr\chi=0$}. Now  we use the automorphisms generated by
boosts, space and time translations to simplify the Jacobi identities. For
example in the  cases (a),(d) and~(e) we may use the boost to put
$\vec{\gamma}=0$\ from the very beginning. On the other hand in the case
(b) by solving Jacobi identity for $\tilde{H},\tilde{P_i},\tilde{P_k}$\ 
we find $\vec{\gamma}\vec{\alpha}=0$\ and again the boost can be used to put
$\vec{\gamma}=0$. Following this way we obtain the eighteen families of solutions 
described in sec.\ref{s3}.

{\Large\bf Acknowledgment}

The authors acknowledge Profs. P.Kosi\'nski  and S.Giller for many helpful discussions.

\end{document}